\documentclass{amsart}

\title{Enumeration of totally positive Grassmann cells}

\author{Lauren K. Williams}

\address{Department of Mathematics, MIT, Cambridge, MA 02139}

\usepackage{pstricks, pst-node, amssymb}

\psset{unit=1pt, arrowsize=4pt, linewidth=.7pt}
\psset{linecolor=blue}
\newgray{grayish}{.90}
\newrgbcolor{embgreen}{0 .5 0}
\def\Le{\hbox{\rotatedown{$\Gamma$}}}
\def\vblack(#1, #2)#3{\cnode*[linecolor=black](#1, #2){3}{#3}}
\def\vwhite(#1,#2)#3{\cnode[linecolor=black,fillcolor=white,fillstyle=solid](#1,#2){3}{#3}}
\countdef\x=23
\countdef\y=24
\countdef\z=25
\countdef\t=26

\def\tbox(#1,#2)#3{
\x=#1 \y=#2 
\multiply\x by 12 
\multiply\y by 12 
\z=\x \t=\y
\advance\z by 12 
\advance\t by 12 
\psline(\x,\y)(\x,\t)(\z,\t)(\z,\y)(\x,\y)
\advance\x by 6
\advance\y by 6 
\rput(\x,\y){{\bf #3}}}

\usepackage{amsmath, amsthm, amssymb, amsbsy}
\usepackage{amsfonts, latexsym, stmaryrd, amscd, xy}
\usepackage[mathscr]{eucal}
\usepackage{epsfig}
\newtheorem{theorem}{Theorem}[section]

\newtheorem{proposition}[theorem]{Proposition}
\newtheorem{lemma}[theorem]{Lemma}

\newtheorem{corollary}[theorem]{Corollary}
\newtheorem{conjecture}[theorem]{Conjecture}
\newtheorem{remark}[theorem]{Remark}

\usepackage{amsfonts}
\usepackage{xy}
\usepackage{amssymb}
\usepackage{amsmath}
\newcommand{\Z}{\mathbb Z}

\newcommand{\R}{\mathbb R}

\DeclareMathOperator{\CB}{CB}

\DeclareMathOperator{\numparts}{numparts}
\DeclareMathOperator{\Arc}{Arc}
\DeclareMathOperator{\modulo}{modulo}
\DeclareMathOperator{\Chord}{Chord}

\DeclareMathOperator{\GL}{GL}
\DeclareMathOperator{\Mat}{Mat}

\newcommand{\thmrefer}[1]{\renewcommand\thetheorem
  {\protect\ref{#1}}\addtocounter{theorem}{-1}}

\xyoption{all}
\CompileMatrices

\begin{document}

\maketitle

\begin{abstract}
Postnikov \cite{Postnikov} has given
a combinatorially explicit 
cell decomposition of  
the totally nonnegative part of a Grassmannian, denoted
$Gr_{k,n}^+$, and showed that this 
set of cells is isomorphic as a graded poset to many other interesting
graded posets.   The main result of our work is an explicit generating
function which enumerates the cells in $Gr_{k,n}^+$ 
 according to their dimension.  As a corollary, we give 
a new proof 
that the Euler characteristic of $Gr_{k,n}^+$ is $1$.
Additionally, we use our result
to produce a  new $q$-analog of the Eulerian numbers,
which interpolates between the Eulerian numbers, the Narayana numbers,
and the binomial coefficients.
\end{abstract}

\section{Introduction}

The classical theory of total positivity concerns matrices 
in which all minors are nonnegative.  
While this theory was pioneered by 
Gantmacher, Krein, and Schoenberg in the 1930s, the past decade has seen
a flurry of research in this area initiated by Lusztig \cite{Lusztig1,
Lusztig2, Lusztig3}. 
Motivated by surprising connections he discovered between 
his theory of canonical bases for quantum groups and the theory of
total 
positivity,
Lusztig 
extended this subject by introducing the totally nonnegative 
variety $G_{\geq 0}$ in an arbitrary reductive group $G$ and
the totally nonnegative part $B_{\geq 0}$ of a real flag variety $B$.
A few years later,
Fomin and Zelevinsky \cite{FZ1} advanced the understanding of 
$G_{\geq 0}$ by studying the 
decomposition of $G$ into double Bruhat cells, and Rietsch \cite{Rietsch}
proved Lusztig's conjectural cell decomposition of $B_{\geq 0}$.
Most recently, Postnikov \cite{Postnikov} investigated the 
combinatorics of the totally nonnegative part of a  
Grassmannian $Gr_{k,n}^+$:  
he established a relationship between 
$Gr_{k,n}^+$ and planar oriented networks, 
producing a combinatorially explicit cell decomposition of 
$Gr_{k,n}^+$.
In this paper we continue Postnikov's study 
of the combinatorics of $Gr_{k,n}^+$:
in particular, we enumerate the cells in the cell decomposition of 
$Gr_{k,n}^+$ according to their dimension.

The {\it totally nonnegative part} of the Grassmannian of $k$-dimensional
subspaces in $\R^n$ is defined to be the quotient 
$Gr_{k,n}^{+} = \GL_k^{+} \backslash \Mat^{+}(k,n)$, where 
$\Mat^{+}(k,n)$ is the space of real
$k \times n$-matrices of rank $k$ with nonnegative maximal minors
and $\GL_k^{+}$ is the group of real matrices with positive determinant.
If we specify which maximal minors are strictly positive and which
are equal to zero, we obtain a cellular decomposition of $Gr_{k,n}^{+}$,
as shown in \cite{Postnikov}.
We refer to the cells in this decomposition as {\it totally positive 
cells}.  The set of totally positive cells naturally has the structure
of a graded poset:
we say that one cell covers another if the closure of the first
cell contains the second, and the rank function is the dimension of 
each cell.

Lusztig \cite{Lusztig1} has proved that the totally nonnegative part of the
(full) flag variety is contractible, which implies the same result 
for any partial flag variety.  (We thank K. Rietsch for pointing this out
to us.)  The topology of the individual cells 
is not well understood, however.  Postnikov \cite{Postnikov} has conjectured
that the closure of each cell in $Gr_{k,n}^+$ is homeomorphic to 
a closed ball. 

In \cite{Postnikov}, Postnikov 
 constructed many different combinatorial objects which 
 are in one-to-one correspondence with the
totally 
positive
Grassmann cells (these objects thereby inherit 
the structure of a graded poset).  
Some of these objects include
decorated permutations, $\Le$-diagrams, positive oriented matroids, and
move-equivalence classes of planar oriented networks. 
Because it is
simple to compute the rank of a particular $\Le$-diagram or decorated
permutation, we will restrict our attention to these two classes
of objects.

The main result of this paper is an explicit formula for the 
{\it rank generating function} $A_{k,n}(q)$ of 
$Gr_{k,n}^+$.  Specifically, 
$A_{k,n}(q)$ is defined to be the
 polynomial in $q$ whose $q^r$ coefficient is the 
number of totally positive cells in $Gr_{k,n}^{+}$ which have dimension $r$. 
As a corollary of our main result,
we give a new proof that
the Euler characteristic of $Gr_{k,n}^+$ is $1$. 
Additionally, using our result and
exploiting the connection between totally positive cells
and permutations,
we compute generating functions
which enumerate (regular) permutations according to two statistics.
This leads to a new $q$-analog of the 
Eulerian numbers that has many interesting combinatorial properties.
For example, when we evaluate this $q$-analog at 
$q=1, 0,$ and $-1$, 
we obtain the Eulerian numbers,
the Narayana numbers, and the binomial coefficients.  Finally, the 
connection with the Narayana numbers suggests a way of incorporating
noncrossing partitions into a larger family of ``crossing" partitions.

Let us fix some notation.
Throughout this paper we use 
$[i]$ to denote the $q$-analog of $i$, that is, 
$[i] = 1+q+ \dots + q^{i-1}$.  (We will sometimes use $[n]$ to refer
to the set $\{1, \dots , n\}$, but the context should make 
our meaning clear.)  Additionally, 
$[i]! := \prod_{k=1}^i [k]$ and 
$\left[ \begin{matrix} i \\ j \end{matrix} \right] := 
\frac{ [i]!}{[j]! [i-j]!}$ are the $q$-analogs of $i!$ and  
 ${i \choose j}$, respectively.

\textsc{Acknowledgments:}
I thank Alex Postnikov for suggesting this problem to me,
and for many helpful 
discussions.  I am indebted to my advisor
Richard Stanley for his invaluable 
advice and constant encouragement. And I 
thank Ira Gessel, Christian Krattenthaler, and 
Konni Rietsch for their very useful comments.
\section{$\Le$-Diagrams}

\begin{figure}[h]
\centerline{\epsfig{figure=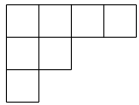}}
\caption{A Young diagram of shape $(4,2,1)$}
\label{YoungDiagram}
\end{figure}

A {\it partition} $\lambda = (\lambda_1, \dots , \lambda_k)$
is a weakly decreasing sequence of nonnegative numbers.  
For a partition $\lambda$, where $\sum \lambda_i = n$, 
the {\it Young diagram} $Y_{\lambda}$ of 
shape $\lambda$ is 
a left-justified diagram of $n$ boxes, with $\lambda_i$ boxes
in the $i$th row.
Figure \ref{YoungDiagram} shows 
a Young diagram of shape $(4,2,1)$.

Fix $k$ and $n$.  Then a {\it $\Le$-diagram} 
$(\lambda, D)_{k,n}$
is a partition $\lambda$ contained in a $k \times (n-k)$ rectangle
(which we will denote by $(n-k)^k$),
together with a filling $D: Y_{\lambda} \to \{0,1\}$ which has the 
{\it $\Le$-property}:
there is no $0$ which has a $1$ above it and a $1$ to its
left.  (Here, ``above" means above and in the same column, and 
``to its left" means to the left and in the same row.)
In Figure \ref{LeDiagram} we give 
an example of a $\Le$-diagram.
\footnote{The symbol $\Le$ is meant to remind the reader of the 
shape of the forbidden pattern, and should be pronounced as
[{le}], because of its relationship to the letter $L$.  See 
\cite{Postnikov} for some interesting numerological remarks on
this symbol.}

\begin{figure}[h]
\pspicture(-120,0)(230,82)

\rput(190,36)
{$\begin{array}{l}
k=6,\ n=17\\
\lambda=(10,9,9,8,5,2)
\end{array}$}

\rput(-10,36){$k$}
\rput(60,82){$n-k$}
\psline[linecolor=black,linewidth=0.5pt]{-}(0,0)(132,0)(132,72)(0,72)(0,0)
\tbox(0,0){1}
\tbox(1,0){1}

\tbox(0,1){0}
\tbox(1,1){0}
\tbox(2,1){0}
\tbox(3,1){1}
\tbox(4,1){1}

\tbox(0,2){0}
\tbox(1,2){0}
\tbox(2,2){0}
\tbox(3,2){0}
\tbox(4,2){0}
\tbox(5,2){0}
\tbox(6,2){1}
\tbox(7,2){1}

\tbox(0,3){0}
\tbox(1,3){0}
\tbox(2,3){0}
\tbox(3,3){0}
\tbox(4,3){0}
\tbox(5,3){0}
\tbox(6,3){0}
\tbox(7,3){0}
\tbox(8,3){0}

\tbox(0,4){1}
\tbox(1,4){1}
\tbox(2,4){1}
\tbox(3,4){1}
\tbox(4,4){0}
\tbox(5,4){1}
\tbox(6,4){1}
\tbox(7,4){1}
\tbox(8,4){1}

\tbox(0,5){0}
\tbox(1,5){1}
\tbox(2,5){1}
\tbox(3,5){0}
\tbox(4,5){0}
\tbox(5,5){1}
\tbox(6,5){0}
\tbox(7,5){1}
\tbox(8,5){0}
\tbox(9,5){1}

\endpspicture

\caption{A $\Le$-diagram $(\lambda, D)_{k,n}$}
\label{LeDiagram}
\end{figure}

We define the rank of 
$(\lambda, D)_{k,n}$ to be the number of 
$1$'s in the filling $D$.
Postnikov proved that there is a one-to-one correspondence between
$\Le$-diagrams $(\lambda, D)$ contained in $(n-k)^k$, 
and totally 
positive cells in $Gr_{k,n}^{+}$, such that
the dimension of a totally positive
cell is equal to the rank of the corresponding $\Le$-diagram.
He proved this by providing a modified 
Gram-Schmidt algorithm $A$, which has the property that it
maps a real $k \times n$ matrix of rank $k$ with nonnegative maximal 
minors to another matrix whose entries are all positive or $0$, which 
has the $\Le$-property.  In brief, the bijection between
totally positive cells and $\Le$-diagrams maps a matrix 
$M$ (representing some totally positive cell) 
to a $\Le$-diagram whose $1$'s represent the positive entries of  
$A(M)$.

Because of this correspondence,
in order to compute $A_{k,n}(q)$, we need to enumerate 
$\Le$-diagrams contained in $(n-k)^k$ according to their number 
of $1$'s.

\section{Decorated Permutations and the  Cyclic Bruhat Order} \label{Bruhat}

The poset of decorated permutations (also called the cyclic Bruhat 
order) was introduced by Postnikov in \cite{Postnikov}. 
A {\it decorated permutation} $\tilde{\pi} = (\pi, d)$ is 
a permutation $\pi$ in the symmetric group
$S_n$ together with a coloring (decoration)
$d$ of its fixed points $\pi (i)=i$ by two colors.
Usually we refer to these two colors as ``clockwise" and 
``counterclockwise," for reasons which the next paragraph will make 
clear.

We represent a decorated permutation $\tilde{\pi}=(\pi,D)$,
where $\pi \in S_n$, by its 
{\it chord diagram}, constructed as follows.  Put $n$ equally
spaced points around a circle, and label these points from $1$ to
$n$ in clockwise order.
If $\pi (i) = j$ then this is represented as a directed arrow, or 
chord,  
from $i$ to $j$.  If $\pi (i) = i$ then
we draw a chord from $i$ to $i$ (i.e. a loop), and orient it either 
clockwise or counterclockwise, according to $d$.
We refer to the chord which begins at position $i$ as 
$\Chord(i)$, and we use $ij$ to denote the directed chord from $i$ to $j$.
Also, if $i,j \in \{1, \dots , n\}$, we use $\Arc(i,j)$ to denote the 
set of points that we would encounter if we were to travel 
clockwise from $i$ to $j$, including $i$ and $j$.

\begin{figure}[h]
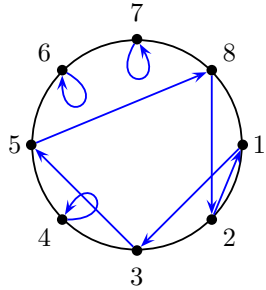

\begin{center}
\pspicture(-60, -60)(60,60)
\pscircle[linecolor=black](0,0){40}
\cnode*[linewidth=0, linecolor=black](40,0){2}{1}
\cnode*[linewidth=0, linecolor=black](28.28,-28.28){2}{2}
\cnode*[linewidth=0, linecolor=black](0, -40){2}{3}
\cnode*[linewidth=0, linecolor=black](-28.28,-28.28){2}{4}
\cnode*[linewidth=0, linecolor=black](-40,0){2}{5}
\cnode*[linewidth=0, linecolor=black](-28.28,28.28){2}{6}
\cnode*[linewidth=0, linecolor=black](0,40){2}{7}
\cnode*[linewidth=0, linecolor=black](28.28,28.28){2}{8}
\rput(46.30,0){$1$}
\rput(35,-35){$2$}
\rput(0,-50){$3$}
\rput(-35,-35 ){$4$}
\rput(-46.30,0){$5$}
\rput(-35,35){$6$}
\rput(0, 50 ){$7$}
\rput(35,35){$8$}
\ncline{->}{1}{3}
\ncline{->}{2}{1}
\ncline{->}{3}{5}
\ncline{->}{5}{8}
\ncline{->}{8}{2}
\nccurve[angleA=-30,angleB=-90,ncurv=20]{->}{6}{6}
\nccurve[angleA=0,angleB=60,ncurv=20]{->}{4}{4}
\nccurve[angleA=-120,angleB=-60,ncurv=20]{->}{7}{7}
\endpspicture
\end{center}
\caption{A chord diagram for a decorated permutation}
\label{chorddiagram}
\end{figure}

For example, the decorated permutation 
$(3,1,5,4,8,6, 7, 2)$ (written in list notation)
with the fixed points $4$, $6$, and $7$ colored in
counterclockwise, clockwise, and counterclockwise, 
respectively, is represented by 
the chord diagram in Figure \ref{chorddiagram}.

The symmetric group $S_n$ acts on the permutations in $S_n$ by 
conjugation.  This action naturally extends to an action of $S_n$
on decorated permutations, if we specify that the action of $S_n$
sends a clockwise (respectively, counterclockwise) fixed point to
a clockwise (respectively, counterclockwise) fixed point.

We say that a pair of chords in a chord diagram
forms a {\it crossing} if they intersect inside the circle or on its
boundary.

\begin{figure}[h]
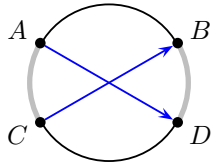

\begin{center}
\pspicture(-50,-40)(50,40)
\pscircle[linecolor=black](0,0){30}
\psarc[linecolor=lightgray,linewidth=2pt](0,0){30}{150}{-150}
\psarc[linecolor=lightgray,linewidth=2pt](0,0){30}{-30}{30}
\cnode*[linewidth=0,linecolor=black](25.98,15){2}{B1}
\cnode*[linewidth=0,linecolor=black](25.98,-15){2}{B2}
\cnode*[linewidth=0,linecolor=black](-25.98,-15){2}{B4}
\cnode*[linewidth=0,linecolor=black](-25.98,15){2}{B5}
\ncline{->}{B5}{B2}
\ncline{->}{B4}{B1}
\rput(-34.64,20){$A$}
\rput(34.64,20){$B$}
\rput(-34.64,-20){$C$}
\rput(34.64,-20){$D$}
\endpspicture
\end{center}
\caption{A crossing}
\label{crossing}
\end{figure}

Every crossing looks like Figure \ref{crossing}, 
where the point $A$ may coincide with the point $B$, and
the point $C$ may coincide with the point $D$.
A crossing is called a {\it simple crossing\/} if there are no
other chords that go from $\Arc(C,A)$ to $\Arc(B,D)$.
Say that two chords are {\it crossing\/} if they form a crossing.

\begin{figure}[h]
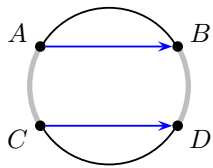

\begin{center}
\pspicture(-50,-40)(50,40)
\pscircle[linecolor=black](0,0){30}
\psarc[linecolor=lightgray,linewidth=2pt](0,0){30}{150}{-150}
\psarc[linecolor=lightgray,linewidth=2pt](0,0){30}{-30}{30}
\cnode*[linewidth=0,linecolor=black](25.98,15){2}{B1}
\cnode*[linewidth=0,linecolor=black](25.98,-15){2}{B2}
\cnode*[linewidth=0,linecolor=black](-25.98,-15){2}{B4}
\cnode*[linewidth=0,linecolor=black](-25.98,15){2}{B5}
\ncline{->}{B5}{B1}
\ncline{->}{B4}{B2}
\rput(-34.64,20){$A$}
\rput(34.64,20){$B$}
\rput(-34.64,-20){$C$}
\rput(34.64,-20){$D$}
\endpspicture
\end{center}
\caption{An alignment}
\label{alignment}
\end{figure}

Let us also say that a pair of chords in a chord diagram forms an
{\it alignment\/} if they are not crossing and they are
relatively located as in Figure \ref{alignment}. 
Here,
again, the point $A$ may coincide with the point $B$,
and the point $C$ may coincide with the point $D$.
If $A$ coincides with $B$ then the chord from $A$ to $B$
should be a counterclockwise
loop in order to be considered an alignment with $\Chord(C)$. 
(Imagine what would happen if we had a piece of string pointing 
from $A$ to $B$, and then we moved the point $B$ to $A$).
And if $C$ coincides with $D$ then the chord from 
$C$ to $D$  should be a
clockwise loop in order to be considered an alignment with 
$\Chord(A)$.  As before, an alignment is a {\it simple alignment\/}
if there are no other chords that go from $\Arc(C,A)$ to
$\Arc(B,D)$.
We say that two chords are {\it aligned\/} if
they form an alignment.

\begin{figure}[h]
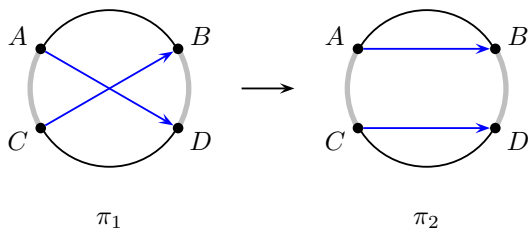

\begin{center}
\pspicture(-50,-60)(170,40)
\pscircle[linecolor=black](0,0){30}
\psarc[linecolor=lightgray,linewidth=2pt](0,0){30}{150}{-150}
\psarc[linecolor=lightgray,linewidth=2pt](0,0){30}{-30}{30}
\rput(-34.64,20){$A$}
\rput(34.64,20){$B$}
\rput(-34.64,-20){$C$}
\rput(34.64,-20){$D$}
\cnode*[linewidth=0,linecolor=black](25.98,15){2}{B1}
\cnode*[linewidth=0,linecolor=black](25.98,-15){2}{B2}
\cnode*[linewidth=0,linecolor=black](-25.98,-15){2}{B4}
\cnode*[linewidth=0,linecolor=black](-25.98,15){2}{B5}
\ncline{->}{B5}{B2}
\ncline{->}{B4}{B1}
\rput(0,-50){$\pi_1$}
\psline[linecolor=black]{->}(50,0)(70,0)
\pscircle[linecolor=black](120,0){30}
\psarc[linecolor=lightgray,linewidth=2pt](120,0){30}{150}{-150}
\psarc[linecolor=lightgray,linewidth=2pt](120,0){30}{-30}{30}
\rput(85.36,20){$A$}
\rput(154.64,20){$B$}
\rput(85.36,-20){$C$}
\rput(154.64,-20){$D$}
\cnode*[linewidth=0,linecolor=black](145.98,15){2}{B21}
\cnode*[linewidth=0,linecolor=black](145.98,-15){2}{B22}
\cnode*[linewidth=0,linecolor=black](94.02,-15){2}{B24}
\cnode*[linewidth=0,linecolor=black](94.02,15){2}{B25}
\ncline{->}{B25}{B21}
\ncline{->}{B24}{B22}
\rput(120,-50){$\pi_2$}
\endpspicture
\end{center}
\caption{Covering relation}
\label{Cover}
\end{figure}

\begin{figure}[h]
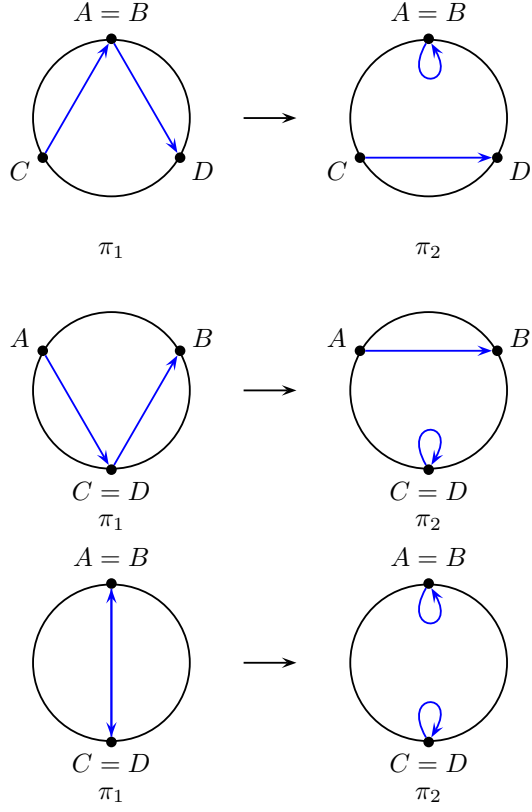

\begin{center}
\pspicture(-50,-60)(170,40)
\pscircle[linecolor=black](0,0){30}
\rput(0,40){$A=B$}
\rput(-34.64,-20){$C$}
\rput(34.64,-20){$D$}
\cnode*[linewidth=0,linecolor=black](0,30){2}{B1}
\cnode*[linewidth=0,linecolor=black](25.98,-15){2}{B2}
\cnode*[linewidth=0,linecolor=black](-25.98,-15){2}{B4}
\ncline{->}{B1}{B2}
\ncline{->}{B4}{B1}
\rput(0,-50){$\pi_1$}
\psline[linecolor=black]{->}(50,0)(70,0)
\pscircle[linecolor=black](120,0){30}
\rput(120,40){$A=B$}
\rput(85.36,-20){$C$}
\rput(154.64,-20){$D$}
\cnode*[linewidth=0,linecolor=black](120,30){2}{B21}
\cnode*[linewidth=0,linecolor=black](145.98,-15){2}{B22}
\cnode*[linewidth=0,linecolor=black](94.02,-15){2}{B24}
\ncline{->}{B24}{B22}
\nccurve[angleA=-120,angleB=-60,ncurv=20]{->}{B21}{B21}
\rput(120,-50){$\pi_2$}
\endpspicture
\end{center}

\begin{center}
\pspicture(-50,-60)(170,40)
\pscircle[linecolor=black](0,0){30}
\rput(0,-38){$C=D$}
\rput(-34.64,20){$A$}
\rput(34.64,20){$B$}
\cnode*[linewidth=0,linecolor=black](0,-30){2}{B1}
\cnode*[linewidth=0,linecolor=black](25.98,15){2}{B2}
\cnode*[linewidth=0,linecolor=black](-25.98,15){2}{B4}
\ncline{->}{B1}{B2}
\ncline{->}{B4}{B1}
\rput(0,-50){$\pi_1$}
\psline[linecolor=black]{->}(50,0)(70,0)
\pscircle[linecolor=black](120,0){30}
\rput(120,-38){$C=D$}
\rput(85.36,20){$A$}
\rput(154.64,20){$B$}
\cnode*[linewidth=0,linecolor=black](120,-30){2}{B21}
\cnode*[linewidth=0,linecolor=black](145.98,15){2}{B22}
\cnode*[linewidth=0,linecolor=black](94.02,15){2}{B24}
\ncline{->}{B24}{B22}
\nccurve[angleA=120,angleB=60,ncurv=20]{->}{B21}{B21}
\rput(120,-50){$\pi_2$}
\endpspicture
\end{center}

\begin{center}
\pspicture(-50,-60)(170,40)
\pscircle[linecolor=black](0,0){30}
\rput(0,40){$A=B$}
\rput(0,-38){$C=D$}
\cnode*[linewidth=0,linecolor=black](0,30){2}{B1}
\cnode*[linewidth=0,linecolor=black](0,-30){2}{B2}
\ncline{->}{B1}{B2}
\ncline{->}{B2}{B1}
\rput(0,-50){$\pi_1$}
\psline[linecolor=black]{->}(50,0)(70,0)
\pscircle[linecolor=black](120,0){30}
\rput(120,40){$A=B$}
\rput(120,-38){$C=D$}
\cnode*[linewidth=0,linecolor=black](120,30){2}{B21}
\cnode*[linewidth=0,linecolor=black](120,-30){2}{B22}
\nccurve[angleA=-120,angleB=-60,ncurv=20]{->}{B21}{B21}
\nccurve[angleA=120,angleB=60,ncurv=20]{->}{B22}{B22}
\rput(120,-50){$\pi_2$}
\endpspicture
\end{center}
\caption{Degenerate covering relations}
\label{Degenerate}
\end{figure}

We now define a partial order on the set of decorated
permutations.
For two decorated permutations $\pi_1$ and $\pi_2$
of the same size $n$, we say that $\pi_1$  {\it covers\/}
$\pi_2$, and write $\pi_1\to\pi_2$, if
the chord diagram of $\pi_1$ contains a pair of chords
that forms a simple crossing 
and the chord diagram of $\pi_2$ is obtained by
changing them to the pair of chords that forms a simple alignment
(see Figure \ref{Cover}).
If the points $A$ and $B$ happen to coincide then the
chord from $A$ to $B$ in the chord diagram of $\pi_2$ degenerates to a
counterclockwise loop.  And if the points $C$ and $D$ coincide then
the chord from $C$ to $D$ in the chord diagram of $\pi_2$ becomes a
clockwise loop.  These degenerate situations are illustrated in 
Figure \ref{Degenerate}. 


Let us define two statistics $A$ and $K$ on decorated permutations.
For a decorated permutation $\pi$, the numbers
$A(\pi)$ and $K(\pi)$ are given by
$$
\begin{array}{l}
A(\pi) = \#\{\textrm{pairs of chords forming an alignment}\},
\\[.1in]
K(\pi) = \#\{i\mid \pi(i)>i\} +
\#\{\textrm{counterclockwise loops}\}.
\end{array}
$$
In our previous example $\pi=(3,1,5,4,8,6,7,2)$
we have $A=11$ and $K=5$.
The $11$  alignments in $\pi$ are
$(13, 66)$, $(21,35)$, $(21,58)$, $(21,44)$,
$(21,77)$, $(35,44)$,
$(35,66)$, $(44,66)$, $(58,77)$, $(66,77)$, $(66,82)$.

\begin{lemma} \cite{Postnikov}
If $\pi_1$ covers $\pi_2$ then $A(\pi_1)=A(\pi_2)-1$
and $K(\pi_1)=K(\pi_2)$.
\label{lem:KA}
\end{lemma}

Note that if $\pi_1$ covers $\pi_2$ then the number of
crossings in $\pi_1$ is greater then the
number of crossings in $\pi_2$.  But the difference of these
numbers is not always $1$.

Lemma~\ref{lem:KA} implies that the transitive closure of the
covering relation ``$\to$'' has the structure of a partially ordered
set and this partially ordered set decomposes into $n+1$ incomparable
components.  
For $0 \leq k\leq n$, we define the {\it cyclic Bruhat order\/}
$\CB_{kn}$ as the set of all decorated permutations $\pi$ of size $n$
such that $K(\pi)=k$  with the partial order relation obtained by
the transitive closure of the covering relation ``$\to$''.
By Lemma~\ref{lem:KA} the function $A$ is the corank function for the cyclic
Bruhat order $\CB_{kn}$.

The definitions of the covering relation and of the statistic
$A$ will not change if we rotate a chord diagram.  The definition
of $K$ depends on the order of the boundary points $1,\dots,n$,
but it is not hard to see that
the statistic $K$ is invariant under the cyclic shift
$\textrm{conj}_\sigma$ for the long cycle $\sigma=(1,2,\dots,n)$.
Thus the order $\CB_{kn}$
is invariant under the action of the cyclic group $\Z/n\Z$
on decorated permutations.

In \cite{Postnikov}, Postnikov proved that the number of totally 
positive cells in $Gr_{k,n}^{+}$ of dimension $r$ is equal to the number
of decorated permutations in $\CB_{kn}$ of rank $r$. 
Thus, $A_{k,n}(1)$ is the cardinality of $\CB_{kn}$, and 
the coefficient of $q^{k(n-k)-\ell}$ in $A_{k,n}(q)$
is the number of decorated permutations in $\CB_{kn}$ 
with $\ell$ alignments.

\section{The Rank Generating Function of $Gr_{k,n}^+$}

Recall that the coefficient of $q^r$ in $A_{k,n}(q)$ is the 
number of cells of dimension $r$ in the 
cellular decomposition of $Gr_{k,n}^+$.  
In this section we use the $\Le$-diagrams to 
find an explicit expression for $A_{k,n}(q)$.  
Additionally, we will  find explicit expressions for the 
generating functions 
$A_{k}(q,x):= \sum_n A_{k,n}(q) x^n$ and 
$A(q,x,y):= \sum_{k\geq 1} \sum_n A_{k,n} (q) x^n y^k$.  
Our main theorem is the following:

\begin{theorem} \label{MainTheorem}
\begin{align*}
A(q,x,y) & = \frac{-y}{q(1-x)} + 
      \sum_{i \geq 1} \frac{y^i (q^{2i+1} -y)}
     {q^{i^2+i+1} (q^i - q^i [i+1] x +[i]xy)} \\
A_k(q,x) &= \sum_{i=0}^{k-1} (-1)^{i+k} 
	  \frac{x^{k-i-1}[i]^{k-i-1}}{q^{ki+i+1}(1-[i+1]x)^{k-i}} +   
       \sum_{i=0}^k (-1)^{i+k}   
	  \frac{x^{k-i}[i]^{k-i}}{q^{ki}(1-[i+1]x)^{k-i+1}} \\
A_{k,n}(q) &= q^{-k^2} \sum_{i=0}^{k-1} (-1)^i {n \choose i}
  (q^{ki} [k-i]^i  [k-i+1]^{n-i} - q^{(k+1)i} [k-i-1]^i [k-i]^{n-i})\\
&= \sum_{i=0}^{k-1} {n \choose i} q^{-(k-i)^2} ([i-k]^i [k-i+1]^{n-i} - 
			   [i-k+1]^i [k-i]^{n-i}).  \\
\end{align*}
\end{theorem}

Note that it is not obvious that $A_{k,n}(q)$ is either polynomial or
nonnegative.

Since the  expressions for $A_k (q,x)$ and 
$A_{k,n}(q)$ follow easily from the
formula for $A(q,x,y)$,  we will concentrate on proving
the formula for $A(q,x,y)$.

\begin{figure}[h]
\centerline{\epsfig{figure=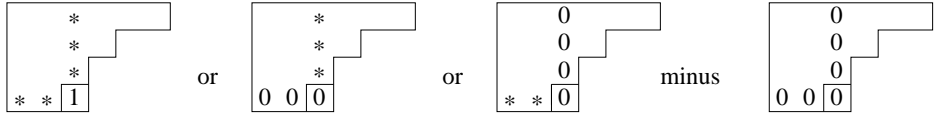}}
\caption{Recurrence for $F_{\lambda}(q)$}
\label{Recurrence}
\end{figure}

Fix a partition $\lambda = (\lambda_1, \dots , \lambda_k)$.  Let
$F_{\lambda} (q)$ be the polynomial in $q$ such that the coefficient of 
$q^r$ is the number of $\Le$-fillings of the Young diagram
$Y_{\lambda}$ which contain $r$ $1$'s.  
As Figure \ref{Recurrence}  illustrates, there is a simple
recurrence for $F_{\lambda}(q)$.   

Explicitly, any 
$\Le$-filling of $\lambda$ is obtained in one of the following ways:
adding a $1$ to the last row of a $\Le$-filling of 
$(\lambda_1, \lambda_2, \dots , \lambda_{k-1}, \lambda_{k} - 1)$;
adding a row containing $\lambda_k$ $0$'s to a $\Le$-filling of 
$(\lambda_1, \dots , \lambda_{k-1})$; or
inserting an all-zero column after the $(\lambda_k - 1)$st column of 
a $\Le$-filling of $(\lambda_1 - 1, \lambda_2 - 1, \dots , \lambda_k - 1)$.
Note, however, that the second and third cases are not exclusive, so 
that our resulting recurrence must subtract off a term corresponding
to their overlap.

\begin{remark} \label{recurrence}
\begin{equation*}
F_\lambda (q) = q F_{(\lambda_1,  \dots , \lambda_{k-1}, \lambda_k -1)}(q) +
    F_{(\lambda_1,  \dots , \lambda_{k-1})}(q) + 
    F_{(\lambda_1 - 1,  \dots , \lambda_k - 1)}(q) - 
    F_{(\lambda_1 - 1,  \dots , \lambda_{k-1} - 1)}(q).
\end{equation*} 
\end{remark}

From the definition, or using the recurrence, it is easy to compute
the first few formulas.  Here are 
$F_{(\lambda_1)}(q)$ and $F_{(\lambda_1, \lambda_2)}(q)$.

\begin{proposition}
\begin{align*} 
F_{(\lambda_1)}(q) &= [2]^{\lambda_1} \\
F_{(\lambda_1,\lambda_2)}(q) &= -q^{-1}[2]^{\lambda_1} +
  q^{-1} [2]^{\lambda_1 - \lambda_2 +1} [3]^{\lambda_2}.
\end{align*}
\end{proposition}

In general, we have the following formula.

\begin{theorem} \label{tool}
Fix $\lambda = (\lambda_1, \lambda_2, \dots , \lambda_k)$.  Then
\begin{equation*}
F_\lambda(q) = 
\sum_{i=1}^k \sum_{1=t_1<\dots < t_i \leq k}
  M(t_1, \dots , t_i : k) 
     [i+1]^{\lambda_{t_i}}
     \prod_{j=2}^{i} [j]^{\lambda_{t_{j-1}}-\lambda_{t_{j}} +1},
\end{equation*}
where $M(t_1, \dots , t_i : k) = (-1)^{k+i} q^{-ik+\sum_{j=1}^i t_j}
       [i]^{k-t_i} 
       \prod_{j=1}^{i-1} [j]^{t_{j+1}-t_j-1}.$
\end{theorem}

Before beginning the proof of the theorem, we state two lemmas which
follow immediately from the formula for $M(t_1, \dots , t_i : k)$.

\begin{lemma}
$M(t_1, \dots , t_i : k) = (-1)^{k-t_i} q^{-i(k-t_i)} [i]^{k-t_i} M(t_1, \dots
, t_i: t_i)$.
\end{lemma}

\begin{lemma} \label{lemma2}
$M(t_1, \dots , t_i : t_i) = -[i-1]^{-1} M(t_1, \dots , t_{i-1} : t_i)$.
\end{lemma}

\begin{proof}
To prove the theorem, we must show that the expression for $F_{\lambda}(q)$
holds for $\lambda=(\lambda_1)$, and that it 
satisfies the recurrence of Remark \ref{recurrence}.  
Also, we must show that 
$F_{(\lambda_1, \lambda_2, \dots , \lambda_k)}(q) = 
F_{(\lambda_1, \lambda_2, \dots , \lambda_k, 0)}(q)$.

The formula  
$F_{(\lambda_1)}(q)=[2]^{\lambda_1}$ clearly agrees with the expression
in the theorem.  
To show that the recurrence is satisfied, we will fix 
$(t_1, \dots , t_i)$ where $1 = t_1 < \dots < t_i \leq k$, and 
calculate the coefficient of 
$[2]^{\lambda_{t_1} -\lambda_{t_2} +1} [3]^{\lambda_{t_2}-\lambda_{t_3}+1} \dots
  [i+1]^{\lambda_{t_i}}$ in each of the five terms of 
  \ref{recurrence}. We will then show that these coefficients
  satisfy the recurrence.

The coefficient in $F_{(\lambda_1, \dots ,\lambda_k)}(q)$ is 
$M(t_1, \dots , t_i: k)$.

The coefficient in 
$F_{(\lambda_1,\lambda_2, \dots ,\lambda_k -1)}(q)$ is
$M(t_1, \dots , t_i : k)$ if $t_i < k$, 
because the term we are looking at together with its coefficient
do not involve $\lambda_k$.
The coefficient is $[i][i+1]^{-1} M(t_1, \dots , t_i :k)$ if $t_i=k$.
    
The coefficient in     
    $F_{(\lambda_1,\lambda_2, \dots ,\lambda_{k-1})}(q)$
is $M(t_1, \dots , t_i : k-1)$ if $t_i<k$, which is equal to 
$-q^i [i]^{-1} M(t_1, \dots , t_i : k)$.  But if $t_i=k$, no such 
term appears, so the coefficient is $0$.

The coefficient in
$F_{(\lambda_1 - 1,\lambda_2 - 1, \dots ,\lambda_k - 1)}(q)$
is always $M(t_1, \dots , t_i: k) [i+1]^{-1}$.

The coefficient in 
$F_{(\lambda_1 - 1,\lambda_2 - 1, \dots ,\lambda_{k-1} - 1)}(q)$ is
$-q^i [i]^{-1} [i+1]^{-1} M(t_1, \dots , t_i: k)$ if $t_i<k$, and $0$
if $t_i = k$.

Let us abbreviate $M(t_1, \dots , t_i:k)$ by $M$.  
We need to show that the coefficients we have just calculated satisfy
the recurrence of Remark \ref{recurrence}.  
For $t_i < k$, this amounts to showing
that
$M = qM - q^i [i]^{-1}M + M[i+1]^{-1} + q^i [i]^{-1} [i+1]^{-1} M$.
And for $t_i = k$, we must show that 
$M = q[i] [i+1]^{-1} M + M[i+1]^{-1}$.  Both of these are easily seen 
to be true.  Thus, we have shown that our expression for $F_{\lambda}(q)$
satisfies Remark \ref{recurrence}.

Now we will show that 
$F_{(\lambda_1, \lambda_2, \dots , \lambda_{k-1}, 0)}(q) = 
F_{(\lambda_1, \lambda_2, \dots , \lambda_{k-1})}(q)$.
It is sufficient to show that the coefficient of 
$[2]^{\lambda_{t_1} - \lambda_{t_2}+1} [3]^{\lambda_{t_2}-\lambda_{t_3}+1} \dots
  [i+1]^{\lambda_{t_i}}$ in $F_{(\lambda_1, \dots , \lambda_k)}(q)$, plus
  $[i+1]$ times the coefficient of 
$[2]^{\lambda_{t_1} - \lambda_{t_2}+1} [3]^{\lambda_{t_2}-\lambda_{t_3}+1} \dots
  [i+1]^{\lambda_{t_i}-\lambda_k + 1} [i+2]^{\lambda_k}$ 
  in $F_{(\lambda_1, \dots , \lambda_k)}(q)$, is equal to the 
  coefficient of $[2]^{\lambda_{t_1}-\lambda_{t_2}+1} \dots 
  [i+1]^{\lambda_{t_i}}$ in $F_{(\lambda_1, \dots , \lambda_{k-1})}(q)$.

In other words, we need
\begin{equation*}
M(t_1, \dots , t_i: k-1) = M(t_1, \dots , t_i: k)+M(t_1, \dots , t_i, k:k)  
  [i+1]. 
\end{equation*}
  From the formula for $M$, 
 we have $M(t_1, \dots , t_i: k-1) = -q^i [i]^{-1} M(t_1, \dots , t_i: k)$.
And from Lemma \ref{lemma2},  
 $M(t_1, \dots , t_i, k:k) = -[i]^{-1}  
   M(t_1, \dots , t_i: k)$.
The proof follows.
\end{proof}

Recall that $A_{k}(q,x)$ is the polynomial in $q$ and $x$ such that
$[q^r x^n] A_{k}(q,x)$ is equal to the 
number of totally positive cells of dimension $r$ in $Gr_{k,n}^+$.  
This is equal to the       
number of $\Le$-diagrams
$(\lambda, D)_{k,n}$ of rank $r$.  We can compute these numbers by
using Theorem \ref{tool}.

\begin{corollary} \label{A_k(q,x)}

\begin{equation*}
A_k(q,x) = \sum_{i=1}^k \sum_{1=t_1 <  \dots < t_{i+1}= k+1} 
  \frac{(-1)^{k+i} q^{-ik+\sum_{j=1}^i t_j} x^k} { (1-x)}
  \prod_{j=1}^i 
  \left(\frac{[j]}{1-[j+1]x}\right)^{t_{j+1}-t_j}. 
\end{equation*}

\end{corollary}

To compute $A_k(q,x)$, we must sum $F_{(\lambda_1, \dots , 
\lambda_k)}(q) x^n$, as $\lambda$ varies over all partitions which fit into 
a $k \times (n-k)$ rectangle.  To do this, we use the following simple lemmas,
the second of which follows immediately from the first.

\begin{lemma} \label{easy}
\begin{equation*}
\sum_{\lambda_1 = 0}^{\infty} \sum_{\lambda_2 = 0}^{\lambda_1} \dots 
 \sum_{\lambda_d = 0}^{\lambda_{d-1}} x_1 ^{\lambda_1} x_2 ^{\lambda_2} \dots
 x_d ^{\lambda_d} = \frac{1}{(1-x_1)(1-x_1 x_2) \dots (1-x_1 x_2 \dots x_d)}.
\end{equation*}
\end{lemma}

\begin{lemma} \label{main}
Fix a set of positive integers
$t_1 < t_2 < \dots < t_d < n+1$.  Then 
\begin{equation*}
\sum_{n=0}^{\infty} \sum_{\lambda_1 = 0}^n \sum_{\lambda_2 = 0}^{\lambda_1} 
\dots \sum_{\lambda_d = 0}^{\lambda_{d-1}} 
[2]^{\lambda_{t_1}-\lambda_{t_2}}  \dots
   [d]^{\lambda_{t_{d-1}}-\lambda_{t_d}} [d+1]^{\lambda_{t_d}} x^n
   \end{equation*} is equal to 
\begin{equation*}
\frac{1}{(1-x)(1-[2]x)^{t_2-t_1}  \dots
(1-[d]x)^{t_d-t_{d-1}} (1-[d+1]x)^{n+1-t_d}}.
\end{equation*}
\end{lemma}

\begin{proof}
For the proof of the corollary, apply 
 Theorem \ref{tool} and Lemma \ref{main} to the fact that
\begin{equation*}A_k(q,x) = 
\sum_{m=0}^{\infty} \sum_{\lambda_1 = 0}^m \sum_{\lambda_2 = 0}^{\lambda_1} 
\dots \sum_{\lambda_k = 0}^{\lambda_{k-1}} F_{(\lambda_1, \dots ,
\lambda_k)}(q) x^m.
\end{equation*} 
\end{proof}

\begin{corollary}
The Euler characteristic of the totally non-negative part of the 
Grassmannian $Gr_{k,n}^+$ is $1$.
\end{corollary}

\begin{proof}
Recall that the Euler characteristic of a cell complex is 
defined to be $\sum_i (-1)^i f_i$, where $f_i$ is the number of 
cells of dimension $i$.  
So if we set $q=-1$ in Corollary \ref{A_k(q,x)}, we will obtain
a polynomial in 
$x$ such that the coefficient of $x^n$ is the Euler characteristic of 
$Gr_{k,n}^+$.
Notice that $[i]$ is equal to $0$ if $i$ is even, and $1$ if $i$ is odd.
So all terms of $A_k(-1,x)$ vanish except the term for $i=1$, which 
becomes $\frac{x^k}{1-x}=x^k+x^{k+1}+x^{k+2}+\dots$.
\end{proof}

Note that this corollary also follows from Lusztig's result that the 
totally nonnegative part of a real flag variety is contractible.

Now our goal will be to simplify our expressions.  
To do so, it is helpful to work
with the ``master" generating function 
$A(q,x,y):= \sum_{k\geq 1}  A_k(q,x) y^k$.   
As a first step, we compute the following expression for $A(q,x,y)$:

\begin{proposition} \label{A(q,x,y)}
\begin{equation*}
A(q,x,y) = \sum_{i=1}^{\infty} q^i [i]! x^i y^i \prod_{j=0}^i \frac{1}{q^j -
q^j[j+1]x + [j]xy}.
\end{equation*}
\end{proposition}

Note that 
$ \frac{1}{q^j -
q^j[j+1]x + [j]xy}$ is not a well-defined formal power series because it
is not clear how to expand it.  In this paper, for reasons which will
become clear in the following proof, we shall always use 
$ \frac{1}{q^j -
q^j[j+1]x + [j]xy}$ as shorthand for the formal power series whose
expansion is implied by the expression 
\begin{equation*}
\frac{1}{q^j (1-[j+1]x)(1 - \frac{q^{-j}[j]y}{1-[j+1]x})}.
\end{equation*}  See \cite[Example~6.3.4]{Stanley2}  
for remarks on the subtleties of such power series.

\begin{proof}
From Corollary \ref{A_k(q,x)}, we know that $A_k(q,x)$ is equal to 
\begin{equation*}
\frac{(-x)^k}{1-x} \sum_{i=1}^k \sum_{1=t_1 <  \dots < 
    t_{i+1}=k+1} 
  (-1)^{i} q^{-ik+\sum_{j=1}^i t_j} 
\prod_{j=1}^i 
\left(\frac{[j]}{1-[j+1]x}\right)^{t_{j+1}-t_j}. 
\end{equation*}
If we make the substitution $\alpha_j = t_{j+1}-t_j$, we then get
\begin{equation*}
A_k(q,x) 
 = \frac{(-x)^k}{1-x} \sum_{i=1}^k (-1)^i q^i 
       \sum_{\substack{\alpha_j \geq 1 \\ \alpha_1 + \dots +\alpha_i = k}}
       \prod_{j=1}^i
      \left(\frac{[j]}{q^j(1-[j+1]x)}\right)^{\alpha_j}.
\end{equation*}

Now let $f_j(p) = ( \frac{[j]}{q^j (1-[j+1]x)})^p$.  For future use, define
$F_j(y):= \sum_{p\geq 1} f_j(p)y^p,$ which is equal to
$\frac{[j]y}{q^j-q^j [j+1]x-[j]y}$.
We get 
\begin{equation*}
A_k(q,x)= \frac{(-x)^k}{1-x} \sum_{i=1}^k (-1)^i q^i 
       \sum_{\substack{\alpha_j \geq 1 \\ \alpha_1 + \dots +\alpha_i = k}}
        \prod_{j=1}^i f_j (\alpha_j), 
\end{equation*}
 and 
we can now easily compute $A(q,x,y):= \sum_{k\geq 1} A_k(q,x)y^k$.
\begin{align*}
A(q,x,y) &= \frac{1}{1-x} \sum_{k\geq 1} (-x)^k \sum_{i=1}^k (-1)^i q^i 
       \sum_{\substack{\alpha_j \geq 1 \\ \alpha_1 + \dots +\alpha_i = k}}
    \prod_{j=1}^i f_j (\alpha_j) y^{\alpha_j}\\
  &= \frac{1}{1-x} \sum_{i=1}^{\infty} \sum_{k\geq i} 
       \sum_{\substack{\alpha_j \geq 1 \\ \alpha_1 + \dots +\alpha_i = k}}
       (-x)^k (-1)^i q^i 
    \prod_{j=1}^i f_j (\alpha_j) y^{\alpha_j}.
\end{align*}
Actually, we can replace $k \geq i$ above with $k \geq 0$, since
if $k<i$ there will be no set of $\alpha_j$ satisfying the conditions
of the third sum.
So we have 

\begin{align*}
A(q,x,y)  &= \frac{1}{1-x} \sum_{i=1}^{\infty} \sum_{k\geq 0} 
       \sum_{\substack{\alpha_j \geq 1 \\ \alpha_1 + \dots +\alpha_i = k}}
       (-x)^k (-1)^i q^i 
    \prod_{j=1}^i f_j (\alpha_j) y^{\alpha_j}\\
  &= \frac{1}{1-x} \sum_{i=1}^{\infty} (-1)^i q^i \sum_{k\geq 0} 
       \sum_{\substack{\alpha_j \geq 1 \\ \alpha_1 + \dots +\alpha_i = k}}
    \prod_{j=1}^i f_j (\alpha_j) (-xy)^{\alpha_j}.
\end{align*}
Using the definition of $F_j$, we get
\begin{align*}
A(q,x,y)
 &= \frac{1}{1-x} \sum_{i=1}^{\infty} (-1)^i q^i 
    \prod_{j=1}^i F_{j}(-xy)\\
 &= \frac{1}{1-x} \sum_{i=1}^{\infty} (-1)^i q^i 
    \prod_{j=1}^i \frac{-[j]xy}{q^j-q^j [j+1]x + [j]xy} \\
 &= \frac{1}{1-x} \sum_{i=1}^{\infty}  q^i  [i]! x^i y^i
    \prod_{j=1}^i \frac{1}{q^j-q^j [j+1]x + [j]xy} \\
 &= \frac{1}{1-x} \sum_{i=1}^{\infty}  q^i  [i]! x^i y^i
    \prod_{j=1}^i \frac{1}{q^j-q^j [j+1]x + [j]xy} \\
 &= \sum_{i=1}^{\infty}  q^i  [i]! x^i y^i
    \prod_{j=0}^i \frac{1}{q^j-q^j [j+1]x + [j]xy}. 
\end{align*}
\end{proof}

Now we will prove the following identity.  This identity
combined with Proposition \ref{A(q,x,y)} will
complete the proof
of Theorem \ref{MainTheorem}.

\begin{theorem} \label{identity}
\begin{equation*}
\sum_{i=1}^{\infty} q^i [i]! x^i y^i \prod_{j=0}^i \frac{1}{q^j -
q^j[j+1]x + [j]xy}=
\frac{-y}{q(1-x)} +\sum_{i \geq 1} \frac{q^{-i^2-i-1} y^i (q^{2i+1} -y)}
     {q^i - q^i [i+1] x +[i]xy}
\end{equation*}
\end{theorem}

\begin{proof}
Observe that the expression on the right-hand side can be thought of
as a partial fraction expansion in terms of $x$, since all denominators
are distinct, and the numerators are free of $x$.
Also note that the $i$-summand of the left-hand side should be easy 
to express in partial fractions with respect to $x$,
since all factors of the denominator are distinct and the 
$x$-degree of the numerator is smaller than the $x$-degree of the denominator.

Thus, our strategy will be to put the left-hand side into 
partial fractions with respect to $x$, and then show that this 
agrees with the right-hand side.

To this end, 
define $\beta_i(j)$ by the equation
\begin{equation*}
\frac{x^i}{\prod_{j=0}^i 
{q^j - q^j[j+1]x+[j]xy}}
= \sum_{j=0}^i \frac{\beta_i(j)}
{q^j - q^j[j+1]x+[j]xy}.
\end{equation*}
Clearing denominators, we obtain 
\begin{equation} \label{useful}
x^i = 
\sum_{j=0}^i \beta_i(j) 
\prod_{\substack{r=0 \\ r\neq j}}^i 
(q^r -q^r [r+1]x+[r]xy).
\end{equation}

Fix $j$.  Notice that 
$(q^j - q^j[j+1]x+[j]xy)$
vanishes when
$x=\frac{q^j}{q^j[j+1]-[j]y}$,
so substitute
$x=\frac{q^j}{q^j[j+1]-[j]y}$ into (\ref{useful}).
We get
\begin{equation*}
\frac{q^{ij}}{(q^j[j+1]-[j]y)^i } = 
\beta_i(j)  
\prod_{\substack{r=0 \\ r\neq j}}^i 
\frac{q^r (q^j [j+1]-[j]y)+q^j([r]y-q^r[r+1])}{q^j[j+1]-[j]y}.
\end{equation*}
Solving for $\beta_i(j)$ and simplifying, we arrive at 
\begin{equation*}
\beta_i(j) = 
\frac{(-1)^{i+j} q^{\frac{j^2+3j-i^2-3i-2ij}{2}} }
{\displaystyle{[j]![i-j]!\prod_{\substack{r=0 \\ r\neq j}}^i (1-q^{-r-j-1}y)}}.
\end{equation*}

Thus the partial fraction expansion with respect to $x$ of 
the left-hand side of Theorem \ref{identity}
is 
\begin{equation*} 
\sum_{i=1}^{\infty} \sum_{j=0}^i 
\frac{\beta_i(j) q^i [i]! y^i}{q^j - q^j [j+1]x+[j]xy}, 
\end{equation*}
which is 
equal to 
\begin{equation} \label{messy}
\sum_{j=0}^{\infty} 
\frac{ {\displaystyle (-1)^j q^{\frac{j^2+3j}{2}} \sum_{\substack{i \geq j\\ i \neq 0}} 
  \left[ \begin{matrix} i\\ j \end{matrix} \right] 
    q^{-{i+1\choose2}-ij}  (-y)^i  
\prod_{\substack{r=0 \\ r\neq j}}^i 
(1-q^{-r-j-1}y)^{-1}}}     
{q^j - q^j [j+1]x+[j]xy}.
\end{equation}

Now it remains to show that the numerator of 
$(q^j-q^j [j+1]x +[j]xy)$ in (\ref{messy}) is equal to 
the numerator of 
$(q^j-q^j [j+1]x +[j]xy)$ in the right-hand side of Theorem
\ref{identity}.
For $j=0$, we must show that 
\begin{equation}  \label{identity1}
(1-\frac{y}{q}) \sum_{i \geq 1} (-1)^i q^{-{i+1\choose 2}} y^i 
\prod_{r=0}^i (1-q^{-r-1}y)^{-1} = 
  \frac{-y}{q}. 
\end{equation}
And 
for $j>0$, we must show that 
\begin{equation} \label{identity2}
(-1)^j q^{\frac{3j^2+j}{2}} y^{-j} \sum_{i \geq j} 
\left[ \begin{matrix} i \\ j \end{matrix} \right] q^{-{i+1\choose 2}-ij} 
(-y)^{i} 
  \prod_{r=0}^i (1-q^{-r-j-1}y)^{-1} = 1.
\end{equation}

If we  make the substitution 
$q \to q^{-1}$ and $r \to r-1$ into (\ref{identity1})
and then add 
the $i=0$ term to both sides, we obtain 

\begin{equation} \label{identity3}
\sum_{i \geq 0} (-1)^i y^i q^{ i+1 \choose 2} 
\prod_{r=1}^{i+1} \frac{1}{1-q^r y} = 1.
\end{equation}

And if we make the same substitution into (\ref{identity2}), we get
\begin{equation} \label{identity4}
(-1)^j q^{- {j+1\choose 2}} y^{-j} \sum_{i \geq j} (-1)^{i} q^{ {i+1 \choose 2}}
\left[ \begin{matrix} i \\ j \end{matrix} \right] y^{i}
\prod_{r=1}^{i+1} \frac{1}{1-q^{r+j}y} = 1 .
\end{equation}

Since (\ref{identity3}) is a special case of 
(\ref{identity4}), it suffices to prove 
(\ref{identity4}).  We will prove this as a separate lemma below; 
modulo this lemma, we are done.
\end{proof}

\begin{lemma} \label{partition}
\begin{equation*}
(-1)^j q^{-{j+1 \choose 2}} y^{-j} \sum_{i \geq j} (-1)^{i} q^ {i+1 \choose 2} 
\left[ \begin{matrix} i \\ j \end{matrix} \right] y^{i}
\prod_{r=1}^{i+1} \frac{1}{1-q^{r+j}y} = 1.
\end{equation*}
\end{lemma}

\begin{proof}
Christian Krattenthaler has pointed out to us
that this lemma is actually a
special case of the ${}_{1}\phi_1$ summation described in Appendix II.5 of 
\cite{GasperRahman}.  
Here, we give
two additional proofs of this lemma.
The first method  
is to show that the infinite sum actually 
telescopes (we thank Ira Gessel for suggesting this to us).   
The second method is to interpret
the lemma as a statement about partitions, and to
prove it combinatorially.

Let us sketch the first method.  We use induction to show that 
\begin{equation*}
(-1)^j q^{-{j+1 \choose 2}} y^{-j} \sum_{i = j}^{m-1} (-1)^{i} q^ {{i+1 \choose 2}}
\left[ \begin{matrix} i \\ j \end{matrix} \right] y^{i}
\prod_{r=1}^{i+1} \frac{1}{1-q^{r+j}y}
\end{equation*}
is equal to 
\begin{equation*}
1 + \frac{ {\displaystyle (-1)^{m-1} q^{jm+{m+1 \choose 2}} y^m
 \sum_{p=0}^j (-1)^p 
\left[ \begin{matrix} m \\ p \end{matrix} \right] 
q^{{p \choose 2}-pj-pm}
y^{-p} }  }
{\prod_{r=1}^{m} (1-q^{r+j}y)}.
\end{equation*}
Then we take the limit as $m$ goes to $\infty$, obtaining  
the statement of the lemma.  

Now let us give a combinatorial proof of the lemma.  For clarity, 
we prove the $j=0$ case in detail and then explain how to
generalize this proof.

First we claim that 
$(-1)^i y^i q^{i+1 \choose 2} \prod_{r=1}^{i+1} \frac{1}{1-q^ry}$
is a generating function for partitions $\lambda$ with 
$i+1$ parts, all distinct, where the smallest part may be zero.
In this formal power series, the coefficient of 
$y^m q^n$ is equal to the number of such partitions with $m$ columns
and $n$ total boxes.
The generating function is multiplied by $1$ or $-1$, 
according to the parity of the number of rows (including zero).

To prove the claim, note that each term of 
$\prod_{r=1}^{i+1} \frac{1}{1-q^r y}$ corresponds to a 
(normal) partition
where rows have lengths between $1$ and $i+1$, inclusive.  The 
exponent of $y$ enumerates the number of rows and the exponent of 
$q$ enumerates the number of boxes.  Now take the transpose of this
partition, so that it is a partition with exactly $i+1$ rows
(possibly zero).  Now the exponent of $y$ is the length of the longest
row.  Add $i, i-1, \dots , 1$ and $0$ boxes to the first, second, ... ,
and
$(i+1)$st rows, respectively.  Finally we have a partition with
$i+1$ parts, all distinct, where the smallest part may be zero.
Since we've added a total of ${i+1 \choose 2}$ boxes to the original
partition, the generating function for this type of partition is 
$q^{i+1 \choose 2} y^i \prod_{r=1}^{i+1} \frac{1}{1-q^r y}$.
Figure \ref{PartitionDiagram1} illustrates the steps in this paragraph.
In the figure, the rows and columns of the partitions are indicated by 
solid and dashed lines, respectively.

\begin{figure}[h]
\centerline{\epsfig{figure=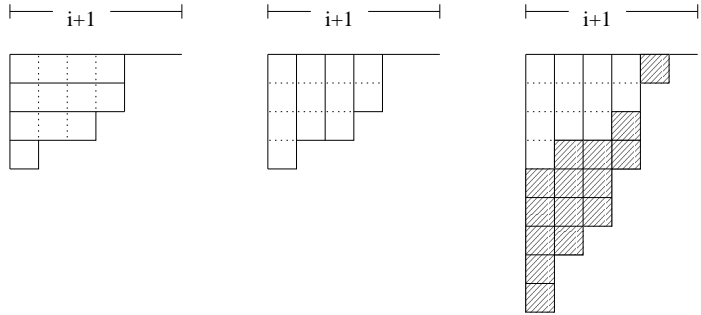}}
\caption{A Combinatorial interpretation for $y^i q^{i+1\choose 2}
  \prod_{r=1}^{i+1} \frac{1}{1-q^r y} $}
\label{PartitionDiagram1} 
\end{figure}

Now we need to find an involution $\phi$ which explains
why all of the terms on the left-hand side of 
(\ref{identity3}) cancel out, except for the $1$.
This involution is very simple: 
if $(\lambda_1, \dots , \lambda_k)$ is a partition 
such that $\lambda_k \neq 0$, then
$\phi(\lambda_1, \dots , \lambda_k) = (\lambda_1, \dots , \lambda_k, 0)$.
And if $\lambda_k = 0$, then
$\phi(\lambda_1, \dots , \lambda_k) = (\lambda_1, \dots , \lambda_{k-1})$.
Clearly both 
$(\lambda_1, \dots , \lambda_k)$ and 
$\phi(\lambda_1, \dots , \lambda_k)$ contribute the 
same powers of $y$ and $q$ to the generating function; the only difference
is the sign.
Only the $0$ partition has no partner under the involution, so all terms
cancel except for $1$.

For the proof of the general case, 
we will show that 
\begin{equation} \label{gf}
q^ {i+1 \choose 2} 
\left[ \begin{matrix} i \\ j \end{matrix} \right] y^{i-j}
\prod_{r=1}^{i+1} \frac{1}{1-q^{r+j}y}
\end{equation}
is a generating function for certain pairs of partitions,
$(\lambda, \hat{\lambda})$.  First, 
note that 
$\prod_{r=1}^{i+1} \frac{1}{1-q^{r+j}y}$ is a generating function
for partitions with rows of lengths $j+1$ through $i+j+1$, inclusive.
It is well-known that 
$\left[ \begin{matrix} i \\ j \end{matrix} \right]$ is a polynomial 
in $q$ whose $q^r$ coefficient 
is the number
of partitions of $r$ which fit inside a $j \times (i-j)$ rectangle.
To account for the 
$\left[ \begin{matrix} i \\ j \end{matrix} \right]$ term in 
(\ref{gf}), let us take 
a partition which fits inside a $j \times (i-j)$ rectangle, and 
place it underneath
a partition with rows of lengths $j+1$ through $i+j+1$, giving us
a partition
with row lengths between $0$ and $i+j+1$, inclusive.  We  
consider this partition to have exactly $i+j+1$ columns, possibly zero.
Finally, to account for the $q^{i+1 \choose 2}$ term in
(\ref{gf}) let us add
$0, 1, \dots , i$ boxes to the last $i+1$ columns of our partition,
so that that the last $i+1$ columns have distinct lengths
(possibly zero).
We now  view the boxes in the 
first $j$ columns of our figure to comprise 
one partition $\lambda$, and the boxes in the 
last $i+1$ columns of our figure
to comprise the transpose of a second partition $\hat{\lambda}$.
Let $\hat{\lambda}_1$ denote the length of the first row of 
$\hat{\lambda}$, and let $r_j(\lambda)$ denote the number of rows
of $\lambda$ which have length $j$.  Then 
the pair $(\lambda, \hat{\lambda})$ satisfies the following conditions:
$\lambda$ has  rows with lengths
between $0$ and $j$, inclusive;
$\hat{\lambda}$ has exactly $i+1$ rows, all distinct,
where the smallest row can have length $0$; and
$r_j(\lambda) +i-j= \hat{\lambda}_1$.
(See Figure 
\ref{PartitionDiagram2} for an illustration of 
$(\lambda, \hat{\lambda})$.)
The term in (\ref{gf}) that corresponds to this 
pair of partitions is $q^{|\lambda|+|\hat{\lambda}|} y^{\numparts (\lambda)}$.

\begin{figure}[h]
\centerline{\epsfig{figure=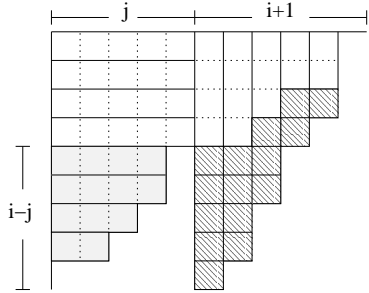}}
\caption{$(\lambda, \hat{\lambda})$, where 
    $\lambda = (5,5,5,5,4,4,3,2,0)$ and $\hat{\lambda}=(9,8,6,4,3,0)$}
\label{PartitionDiagram2} 
\end{figure}
Our involution $\phi$ is a 
simple generalization of the involution we used before.
This time, $\phi$ fixes $\lambda$,  and either adds or subtracts
a trailing zero to $\hat{\lambda}$.

\end{proof}

This completes the proof of Theorem \ref{MainTheorem}.

In Table \ref{A(k,n)}, we have listed some of the values of 
$A_{k,n}(q)$ for small $k$ and $n$.
It is easy to see from the definition of $\Le$-diagrams that
$A_{k,n}(q)=A_{n-k,n}(q)$: one can reflect a $\Le$-diagram
$(\lambda, D)_{k,n}$ of rank $r$ over the main diagonal
to get another $\Le$-diagram $(\lambda^{\prime}, D^{\prime})_{n-k,n}$ of
rank $r$.  Alternatively, one should be able to prove the claim
directly from the expression in Theorem \ref{MainTheorem}, 
using some $q$-analog of 
Abel's identity.
\bigskip
\begin{table}[h]
\begin{tabular}{|p{1.1cm}|p{10.5cm}|}
\hline
$A_{1,1}(q)$  & $1$ \\
$A_{1,2}(q)$  & $q+2$ \\
$A_{1,3}(q)$  & $q^2+3q+3$ \\
$A_{1,4}(q)$  & $q^3+4q^2+6q+4$ \\
\hline
$A_{2,4}(q)$  & $q^4+4q^3+10q^2+12q+6$ \\
$A_{2,5}(q)$  & $q^6+5q^5+15q^4+30q^3+40q^2+30q+10$ \\
$A_{2,6}(q)$  & $q^8+6q^7+21q^6+50q^5+90q^4+120q^3+110q^2+60q+15$ \\
\hline
$A_{3,6}(q)$  & $q^9+6q^8+21q^7+56q^6+114q^5+180q^4+215q^3+180q^2+90q+20$ \\
$A_{3,7}(q)$  &
$q^{12}+7q^{11}+28q^{10}+84q^9+203q^8+406q^7+679q^6+938q^5+1050q^4+
910q^3+560q^2+210q+35$ \\
\hline
\end{tabular}
\bigskip
\caption{$A_{k,n}(q)$}
\label{A(k,n)}
\end{table}

Note that it is possible to see directly from the definition that 
$Gr_{1,n}^+$ is just some deformation of a simplex with $n$ vertices.
This explains the simple form of $A_{1,n}(q)$.

\section{A New $q$-Analog of the Eulerian Numbers}

If $\pi \in S_n$, we say that $\pi$ has a {\it weak excedence} 
at position $i$ if $\pi (i) \geq i$.  
The {\it Eulerian number} $E_{k,n}$ is the number of permutations
in $S_n$ which have $k$ weak excedences.  (One can define the 
Eulerian numbers in terms of other statistics, such as descent, but
this will not concern us here.)

Now that we have computed the rank generating function for 
$\CB_{kn}^+$ (which is the rank generating function for  the poset 
of decorated permutations),
we can use this result to enumerate (regular)
permutations according to two statistics: weak excedences and alignments.
This gives us a new $q$-analog of the Eulerian numbers.

Recall that the statistic $K$ on decorated permutations
was defined as 
$$
K(\pi) = \#\{i\mid \pi(i)>i\} +
\#\{\textrm{counterclockwise loops}\}.
$$
Note that $K$ is related to the notion of weak excedence 
in permutations.  In fact, we can 
extend the definition of weak excedence to decorated permutations by
saying that a decorated permutation has a weak excedence in 
position $i$, if $\pi(i)>i$, or if $\pi(i)=i$ and $d(i)$ is 
counterclockwise.  
This makes sense, since 
the limit of a chord from $1$ to $2$ as $1$ approaches $2$, is a 
counterclockwise loop.  
Then $K(\pi)$ is the number of weak excedences
in $\pi$.

We will call a decorated permutation {\it regular} if all of its
fixed points are oriented counterclockwise.  Thus, a fixed point of
a regular permutation will always be a weak excedence, as it should be. 
Recall that the Eulerian number $E_{k,n}$ is the number of 
permutations of $[n]$ with $k$ weak excedences.
Earlier, we saw
that the coefficient of $q^{k(n-k)-\ell}$ in $A_{k,n}(q)$
is the number of decorated permutations in $\CB_{kn}$ 
with $\ell$ alignments.
By analogy, let $E_{k,n}(q)$ be the 
polynomial in $q$ whose coefficient of 
$q^{k(n-k)-\ell}$ is the number of (regular) permutations with 
$k$ weak excedences and $\ell$ alignments.  
Thus, the family $E_{k,n}(q)$ will be a $q$-analog of 
the Eulerian numbers. 

We can relate decorated permutations to regular permutations
via the following lemma.

\begin{lemma}
$A_{k,n}(q) = \sum_{i=0}^{n} {n \choose i} E_{k,n-i}(q).$
\end{lemma}

\begin{proof}
To prove this lemma we need to figure out how the number of 
alignments changes,
if we start with a regular permutation on $[n-i]$ with
$k$ weak excedences, and then
add $i$ clockwise fixed points.  
Note that adding a clockwise fixed point adds exactly $k$ alignments,
since a clockwise fixed point is aligned with all of the weak excedences.
Since clockwise fixed points are not in alignment with each other, it 
follows that adding $i$ clockwise fixed points adds exactly 
$ik$ alignments.  

This shows that the new number of alignments is equal to $ki$ plus
the old number of alignments, or equivalently, that 
$k(n-i-k)$ minus the old number of alignments is equal to 
$k(n-k)$ minus the new number of alignments.
In other words, the rank of the permutation on $[n-i]$ 
is equal to the rank of the new decorated permutation on $[n]$.
Both permutations have $k$ weak excedences.  
Since there are ${n \choose i}$ ways to pick $i$ entries of 
a permutation on $[n]$ to be designated as clockwise
fixed points, we have that 
$A_{k,n}(q) = \sum_{i=0}^{n} {n \choose i} E_{k,n}(q).$

\end{proof}

Observe that we can invert the formula given in the lemma, deriving the 
following corollary.

\begin{corollary}
\begin{equation*}
E_{k,n}(q) = \sum_{i=0}^n (-1)^{i} {n \choose i} A_{k,n-i}(q).
\end{equation*}
\end{corollary}

Putting this together with Theorem \ref{MainTheorem}, we get the following.

\begin{corollary} \label{qEulerian}
\begin{align*}
E_{k,n}(q) &= q^{n-k^2} \sum_{i=0}^{k-1} {n \choose i} (-1)^i
	      (q^{ki-i}[k-i]^n - q^{ki}[k-i-1]^n)\\ 
   &= q^{n-k^2} \sum_{i=0}^{k-1} (-1)^i [k-i]^n q^{ki-k}
      ( {n \choose i} q^{k-i} + {n \choose i-1}).
\end{align*}
\end{corollary}

Notice that by substituting $q=1$ into the second formula, we get
$$E_{k,n} = 
\sum_{i=0}^k (-1)^i {n+1 \choose i} (k-i)^n,$$ the well-known 
exact formula for the Eulerian numbers.

Now we will investigate the properties of 
$E_{k,n}(q)$.
Actually, since 
$E_{k,n}(q)$ is a multiple of $q^{n-k}$, we first
define 
$\hat{E}_{k,n}(q)$ to be $q^{k-n} E_{k,n}(q)$, and then
work with this renormalized polynomial.
Table \ref{E(k,n)} lists 
 $\hat{E}_{k,n}(q)$ for $n=4, 5, 6, 7$.

\bigskip
\begin{table}[h]
\begin{tabular}{|p{1.1cm}|p{10cm}|}
\hline $\hat{E}_{1,4}(q)$   & $1$ \\
 $\hat{E}_{2,4}(q)$   & $6+4q+q^2$ \\
 $\hat{E}_{3,4}(q)$   & $6+4q+q^2$ \\
 $\hat{E}_{4,4}(q)$   & $1$ \\
\hline $\hat{E}_{1,5}(q)$   & $1$ \\
 $\hat{E}_{2,5}(q)$   & $10+10q+5q^2+q^3$ \\
 $\hat{E}_{3,5}(q)$   & $20+25q+15q^2+5q^3+q^4$ \\
 $\hat{E}_{4,5}(q)$   & $10+10q+5q^2+q^3$ \\
 $\hat{E}_{5,5}(q)$   & $1$ \\
\hline $\hat{E}_{1,6}(q)$   & $1$ \\
 $\hat{E}_{2,6}(q)$   & $15+20q+15q^2+6q^3+q^4$ \\
 $\hat{E}_{3,6}(q)$   & $50+90q+84q^2+50q^3+21q^4+6q^5+q^6$ \\
 $\hat{E}_{4,6}(q)$   & $50+90q+84q^2+50q^3+21q^4+6q^5+q^6$ \\
 $\hat{E}_{5,6}(q)$   & $15+20q+15q^2+6q^3+q^4$ \\
 $\hat{E}_{6,6}(q)$   & $1$ \\
\hline $\hat{E}_{1,7}(q)$ & $1$ \\
 $\hat{E}_{2,7}(q)$   & $21+35q+35q^2+21q^3+7q^4+q^5$ \\
 $\hat{E}_{3,7}(q)$   & $105+245q+308q^2+259q^3+161q^4+77q^5+28q^6+7q^7+q^8$ \\
 $\hat{E}_{4,7}(q)$   & $175+441q+588q^2+532q^3+364q^4+196q^5+84q^6+28q^7+7q^8+q^9$ \\  
 $\hat{E}_{5,7}(q)$   & $105+245q+308q^2+259q^3+161q^4+77q^5+28q^6+7q^7+q^8$ \\
 $\hat{E}_{6,7}(q)$   & $21+35q+35q^2+21q^3+7q^4+q^5$ \\
 $\hat{E}_{7,7}(q)$ & $1$ \\
\hline
\end{tabular}
\bigskip
\caption{$\hat{E}_{k,n}(q)$}
\label{E(k,n)}
\end{table}
\bigskip

We can make a number of observations about these polynomials.  For
example, we can generalize the well-known result 
that $E_{k,n}=E_{n+1-k,n}$, where $E_{k,n}$ is the Eulerian number
corresponding to the number of permutations of $S_n$ with $k$
weak excedences.

\begin{proposition}
$\hat{E}_{k,n}(q) = 
\hat{E}_{n+1-k,n}(q).$ 
\end{proposition}

\begin{proof}
To prove this, we define an alignment-preserving  bijection on the set of 
permutations in $S_n$, which maps permutations with $k$ weak excedences
to permutations with $n+1-k$ weak excedences.
If $\pi =( a_1, a_2, \dots ,a_n)$ is a permutation written in list notation, then
the bijection maps $\pi$ to $(b_1, b_2, \dots ,b_n)$, where  
$b_i = n - a_{n+1-i} \modulo n$.
\end{proof}

The reader will probably have noticed from the table that 
the coefficients of $\hat{E}_{2,n}(q)$ are binomial coefficients.
Indeed, we have the following proposition, which follows from 
Corollary \ref{qEulerian}.

\begin{proposition}
$\hat{E}_{2,n}(q) = \sum_{i=0}^{n-2} {n \choose i+2} q^i$.
\end{proposition}

\begin{proposition}
\cite{Postnikov}
The coefficient of the highest degree term of 
$\hat{E}_{k,n}(q)$
is $1$.  
\end{proposition}
\begin{proof}
This is 
because there is a unique permutation in $S_n$ with $k$ weak excedences
and no alignments, as proved in \cite{Postnikov}.
That unique permutation is $\pi_k : i \mapsto i+k \modulo n$.
\end{proof}

\begin{proposition}  \label{binomial}
$\hat{E}_{k,n}(-1) = \pm {n-1 \choose k-1}$.
\end{proposition}

\begin{proof}
If we substitute $q=-1$ into the first expression for $E_{k,n}(q)$,
we eventually get 
$(-1)^{n+1} \sum_{i=0}^{k-1} {n \choose i} (-1)^i$.  It is known 
(see \cite{Andrews}) that
this expression is equal to ${n-1 \choose k-1}$.

\end{proof}

\begin{proposition} \label{Narayana}
$\hat{E}_{k,n}(q)$ is a polynomial
of degree $(k-1)(n-k)$, and 
$\hat{E}_{k,n}(0)$
is the Narayana number $N_{k,n}= \frac{1}{n} {n \choose k} {n \choose k-1}$.
\end{proposition}

We will prove Proposition \ref{Narayana} in Section \ref{NaryanaConnection}.

\begin{corollary}
$\hat{E}_{k,n}(q)$ interpolates between the Eulerian numbers, the Narayana 
numbers, and the binomial coefficients, at $q=1, 0$, and $-1$, respectively.
\end{corollary}
\begin{proof}
This follows from the
fact that $\hat{E}_{k,n}(q)$ is a $q$-analog of the Eulerian numbers,
together with Propositions \ref{binomial} and \ref{Narayana}.
\end{proof}

Based on experimental evidence, we formulated the following conjecture
about the coefficient of $q$ in 
$\hat{E}_{k,n}(q)$.  However, nice expressions for coefficients of 
other terms
 have eluded us so far. 

\begin{conjecture}
The coefficient of 
$q$ in $\hat{E}_{k,n}(q)$ 
is $ {n \choose k+1} {n \choose k-2}$.
\end{conjecture}

\begin{remark}
The coefficients of $\hat{E}_{k,n}(q)$ appear to be unimodal.  However,
these polynomials do not in general have real zeroes.
\end{remark}

Since it may be helpful to have formulas which 
enumerate permutations by alignments (rather than $k(n-k)$ minus the
number of alignments), we let 
$\widetilde{E}_{k,n}(q)$ be the polynomial in $q$ such that
the coefficient of $q^l$ is the number of permutations on 
$\{1, \dots n \}$ with $k$ weak excedences and $l$ alignments.
Note that by using Corollary \ref{qEulerian} and
performing a transformation  which 
sends $q$ to $q^{-1}$, we get the following
expressions.

\begin{align*}
\widetilde{E}_{k,n}(q) &= \sum_{i=0}^{k-1} {n \choose i} (-1)^i
  q^{i(n-k)} ( q^i [k-i]^n - q^n [k-i-1]^n )\\
  &= 
  \sum_{i=0}^{k-1} (-1)^i [k-i]^n q^{i(n-k)} ( {n \choose i} q^i +
	  { n \choose i-1 } q^k ).
\end{align*}

\section{Connection with Narayana Numbers} \label{NaryanaConnection}

A {\it noncrossing partition} of the set $[n]$ is a partition 
$\pi$ of the set $[n]$ with the property that if $a < b < c < d$
and some block $B$ of $\pi$ contains both $a$ and $c$, while some 
block $B^{\prime}$ of $\pi$ contains both $b$ and $d$, then $B=B^{\prime}$.
Graphically, we can represent a noncrossing partition on a circle
which has $n$ labeled points equally spaced around it.  We represent
each block $B$ as the polygon whose vertices are the elements of $B$.
Then the condition that $\pi$ is noncrossing just means that no two
blocks (polygons) intersect each other.

It is known that the number of 
noncrossing partitions of $[n]$ which have $k$ blocks is equal to 
the Narayana number $N_{k,n}= \frac{1}{n} {n \choose k}{n \choose k-1}$
(see Exercise 68e in \cite{Stanley2}).

To prove the following proposition we will find a bijection between
permutations of $S_n$ with $k$ excedences and the maximal number of alignments,
and noncrossing partitions on $[n]$. 

\begin{proposition}
Fix $k$ and $n$. Then
$(k-1)(n-k)$ is the maximal number of alignments that a permutation
in $S_n$ with $k$ weak excedences can have. 
The number of permutations in $S_n$ with $k$ weak excedences
that achieve the maximal number
of alignments is the Narayana number $N_{k,n} = \frac{1}{n}
{n \choose k}{n \choose k-1}$.
\end{proposition}

\begin{proof}
Recall the bijection between $\Le$-diagrams and decorated permutations.
The $\Le$-diagrams which correspond to regular permutations 
with $k$ weak excedences are the 
$\Le$-diagrams $(\lambda, D)$ contained in a $k$ by $n-k$ rectangle, such
that each column of the rectangle contains at least one $1$.  
The squares of the rectangle which do not contain a $1$ correspond to
alignments, so the maximal number of alignments is $(k-1)(n-k)$.
(It is also straightforward to prove this using decorated permutations.)

In order to prove that the number of permutations which achieve the maximum
number of alignments is $N_{k,n}$, we put these permutations in bijection
with noncrossing partitions of $[n]$ which have $k$ blocks.  

To figure out what the maximal-alignment permutations look like, 
imagine starting from any given permutation and applying the 
covering relations in the cyclic Bruhat order as many times as possible,
such that the result is a regular permutation.  
Note that of the four cases of the covering relation (illustrated in 
section \ref{Bruhat}), we can use only the first and 
second cases.
We cannot use the third and fourth operations because 
these add clockwise fixed points, which are not allowed in regular
permutations.
It is easy to see that  after applying the first two operations as many
times as possible, the resulting permutation will have no crossings
among its chords and all cycles will be directed counterclockwise.

The map from maximal-alignment permutations to noncrossing
partitions is now obvious.  We simply take our permutation and then
erase the directions on the edges.  Since the covering
relations in the cyclic Bruhat order preserve the number of 
weak excedences, and since each counterclockwise cycle in a 
permutation contributes one weak excedence, the resulting 
noncrossing partitions will all have $k$ blocks.  
In Figure \ref{Bijection} 
we show the permutations in $S_4$ which have $2$ weak excedences
and $2$ alignments,
along with the corresponding noncrossing partitions.
\begin{figure}
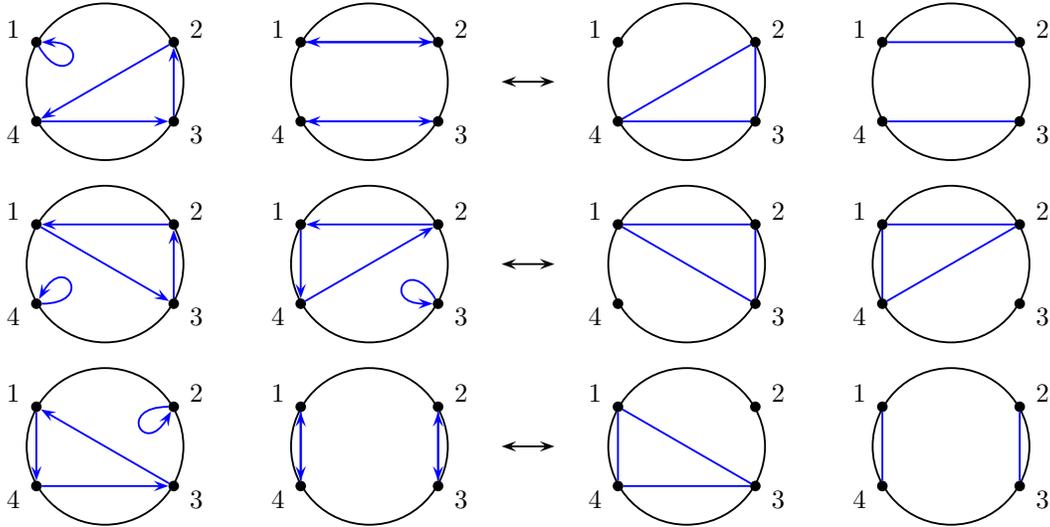

\begin{center}
\pspicture(-120,-10)(340,40)
\pscircle[linecolor=black](-100,0){30}
\rput(-134.64,20){$1$}
\rput(-65.36,20){$2$}
\rput(-134.64,-20){$4$}
\rput(-65.36,-20){$3$}
\cnode*[linewidth=0,linecolor=black](-74.02,15){2}{B2}
\cnode*[linewidth=0,linecolor=black](-74.02,-15){2}{B3}
\cnode*[linewidth=0,linecolor=black](-125.98,-15){2}{B4}
\cnode*[linewidth=0,linecolor=black](-125.98,15){2}{B1}
\ncline{->}{B2}{B4}
\ncline{->}{B3}{B2}
\ncline{->}{B4}{B3}
\nccurve[angleA=-60,angleB=0,ncurv=20]{->}{B1}{B1}
\pscircle[linecolor=black](0,0){30}
\rput(-34.64,20){$1$}
\rput(34.64,20){$2$}
\rput(-34.64,-20){$4$}
\rput(34.64,-20){$3$}
\cnode*[linewidth=0,linecolor=black](25.98,15){2}{B2}
\cnode*[linewidth=0,linecolor=black](25.98,-15){2}{B3}
\cnode*[linewidth=0,linecolor=black](-25.98,-15){2}{B4}
\cnode*[linewidth=0,linecolor=black](-25.98,15){2}{B1}
\ncline{->}{B1}{B2}
\ncline{->}{B2}{B1}
\ncline{->}{B3}{B4}
\ncline{->}{B4}{B3}
\psline[linecolor=black]{<->}(50,0)(70,0)
\pscircle[linecolor=black](120,0){30}
\rput(85.36,20){$1$}
\rput(154.64,20){$2$}
\rput(85.36,-20){$4$}
\rput(154.64,-20){$3$}
\cnode*[linewidth=0,linecolor=black](145.98,15){2}{B2}
\cnode*[linewidth=0,linecolor=black](145.98,-15){2}{B3}
\cnode*[linewidth=0,linecolor=black](94.02,-15){2}{B4}
\cnode*[linewidth=0,linecolor=black](94.02,15){2}{B1}
\ncline{-}{B2}{B4}
\ncline{-}{B3}{B2}
\ncline{-}{B4}{B3}
\pscircle[linecolor=black](220,0){30}
\rput(185.36,20){$1$}
\rput(254.64,20){$2$}
\rput(185.36,-20){$4$}
\rput(254.64,-20){$3$}
\cnode*[linewidth=0,linecolor=black](245.98,15){2}{B2}
\cnode*[linewidth=0,linecolor=black](245.98,-15){2}{B3}
\cnode*[linewidth=0,linecolor=black](194.02,-15){2}{B4}
\cnode*[linewidth=0,linecolor=black](194.02,15){2}{B1}
\ncline{-}{B1}{B2}
\ncline{-}{B3}{B4}
\endpspicture
\end{center}

\begin{center}
\pspicture(-120,-10)(340,40)
\pscircle[linecolor=black](-100,0){30}
\rput(-134.64,20){$1$}
\rput(-65.36,20){$2$}
\rput(-134.64,-20){$4$}
\rput(-65.36,-20){$3$}
\cnode*[linewidth=0,linecolor=black](-74.02,15){2}{B2}
\cnode*[linewidth=0,linecolor=black](-74.02,-15){2}{B3}
\cnode*[linewidth=0,linecolor=black](-125.98,-15){2}{B4}
\cnode*[linewidth=0,linecolor=black](-125.98,15){2}{B1}
\ncline{->}{B2}{B1}
\ncline{->}{B3}{B2}
\ncline{->}{B1}{B3}
\nccurve[angleA=0,angleB=60,ncurv=20]{->}{B4}{B4}
\pscircle[linecolor=black](0,0){30}
\rput(-34.64,20){$1$}
\rput(34.64,20){$2$}
\rput(-34.64,-20){$4$}
\rput(34.64,-20){$3$}
\cnode*[linewidth=0,linecolor=black](25.98,15){2}{B2}
\cnode*[linewidth=0,linecolor=black](25.98,-15){2}{B3}
\cnode*[linewidth=0,linecolor=black](-25.98,-15){2}{B4}
\cnode*[linewidth=0,linecolor=black](-25.98,15){2}{B1}
\ncline{->}{B2}{B1}
\ncline{->}{B4}{B2}
\ncline{->}{B1}{B4}
\nccurve[angleA=120,angleB=180,ncurv=20]{->}{B3}{B3}
\psline[linecolor=black]{<->}(50,0)(70,0)
\pscircle[linecolor=black](120,0){30}
\rput(85.36,20){$1$}
\rput(154.64,20){$2$}
\rput(85.36,-20){$4$}
\rput(154.64,-20){$3$}
\cnode*[linewidth=0,linecolor=black](145.98,15){2}{B2}
\cnode*[linewidth=0,linecolor=black](145.98,-15){2}{B3}
\cnode*[linewidth=0,linecolor=black](94.02,-15){2}{B4}
\cnode*[linewidth=0,linecolor=black](94.02,15){2}{B1}
\ncline{-}{B2}{B1}
\ncline{-}{B3}{B2}
\ncline{-}{B1}{B3}
\pscircle[linecolor=black](220,0){30}
\rput(185.36,20){$1$}
\rput(254.64,20){$2$}
\rput(185.36,-20){$4$}
\rput(254.64,-20){$3$}
\cnode*[linewidth=0,linecolor=black](245.98,15){2}{B2}
\cnode*[linewidth=0,linecolor=black](245.98,-15){2}{B3}
\cnode*[linewidth=0,linecolor=black](194.02,-15){2}{B4}
\cnode*[linewidth=0,linecolor=black](194.02,15){2}{B1}
\ncline{-}{B2}{B1}
\ncline{-}{B4}{B2}
\ncline{-}{B4}{B1}
\endpspicture
\end{center}

\begin{center}
\pspicture(-120,-30)(340,40)
\pscircle[linecolor=black](-100,0){30}
\rput(-134.64,20){$1$}
\rput(-65.36,20){$2$}
\rput(-134.64,-20){$4$}
\rput(-65.36,-20){$3$}
\cnode*[linewidth=0,linecolor=black](-74.02,15){2}{B2}
\cnode*[linewidth=0,linecolor=black](-74.02,-15){2}{B3}
\cnode*[linewidth=0,linecolor=black](-125.98,-15){2}{B4}
\cnode*[linewidth=0,linecolor=black](-125.98,15){2}{B1}
\ncline{->}{B1}{B4}
\ncline{->}{B4}{B3}
\ncline{->}{B3}{B1}
\nccurve[angleA=180,angleB=240,ncurv=20]{->}{B2}{B2}
\pscircle[linecolor=black](0,0){30}
\rput(-34.64,20){$1$}
\rput(34.64,20){$2$}
\rput(-34.64,-20){$4$}
\rput(34.64,-20){$3$}
\cnode*[linewidth=0,linecolor=black](25.98,15){2}{B2}
\cnode*[linewidth=0,linecolor=black](25.98,-15){2}{B3}
\cnode*[linewidth=0,linecolor=black](-25.98,-15){2}{B4}
\cnode*[linewidth=0,linecolor=black](-25.98,15){2}{B1}
\ncline{->}{B1}{B4}
\ncline{->}{B4}{B1}
\ncline{->}{B2}{B3}
\ncline{->}{B3}{B2}
\psline[linecolor=black]{<->}(50,0)(70,0)
\pscircle[linecolor=black](120,0){30}
\rput(85.36,20){$1$}
\rput(154.64,20){$2$}
\rput(85.36,-20){$4$}
\rput(154.64,-20){$3$}
\cnode*[linewidth=0,linecolor=black](145.98,15){2}{B2}
\cnode*[linewidth=0,linecolor=black](145.98,-15){2}{B3}
\cnode*[linewidth=0,linecolor=black](94.02,-15){2}{B4}
\cnode*[linewidth=0,linecolor=black](94.02,15){2}{B1}
\ncline{-}{B3}{B1}
\ncline{-}{B4}{B3}
\ncline{-}{B1}{B4}
\pscircle[linecolor=black](220,0){30}
\rput(185.36,20){$1$}
\rput(254.64,20){$2$}
\rput(185.36,-20){$4$}
\rput(254.64,-20){$3$}
\cnode*[linewidth=0,linecolor=black](245.98,15){2}{B2}
\cnode*[linewidth=0,linecolor=black](245.98,-15){2}{B3}
\cnode*[linewidth=0,linecolor=black](194.02,-15){2}{B4}
\cnode*[linewidth=0,linecolor=black](194.02,15){2}{B1}
\ncline{-}{B4}{B1}
\ncline{-}{B2}{B3}
\endpspicture
\end{center}
\caption{The bijection between maximal-alignment permutations and 
noncrossing partitions}
\label{Bijection}
\end{figure}

Conversely, if we start with a noncrossing partition on $[n]$ which 
has $k$ blocks, and then orient each cycle counterclockwise, then
this gives us a maximal-alignment permutation with $k$ weak
excedences.

The map from maximal-alignment permutations to noncrossing permutations
is obvious.  Note that a maximal-alignment permutation must correspond 
to a noncrossing partition because, if there were a crossing of chords,
we could uncross them to increase the number of alignments
(while preserving the number of excedences).

\end{proof}

\begin{corollary}
The number of permutations in $S_n$ which have the maximal number of 
alignments, given their weak excedences, is 
$C_n = \frac{1}{n} {2n \choose n+1}$, the $n$th Catalan number.
\end{corollary}

\begin{proof}
It is known that $\sum_k N_{k,n} = C_n$.
\end{proof}

\begin{remark}
The bijection between maximal-alignment permutations and noncrossing
partitions is especially interesting because the connection gives 
a way of incorporating 
noncrossing partitions into a larger family of ``crossing" partitions;
this family of crossing partitions is a ranked poset,
graded by alignments.
\end{remark}

\section{Connections with the Permanent}

Let $M_n (x)$ denote the permanent of the $n \times n$ matrix 
\begin{equation*}
\left(
  \begin{matrix}
    1+x & x & x & \dots & x \\
    1 & 1+x & x & \dots & x \\
    1 & 1 & 1+x & \dots & x \\
    \vdots & \vdots & \vdots & & \vdots \\
    1 & 1 & 1 & \dots & 1+x 
  \end{matrix}
  \right).
\end{equation*}

Clearly $[x^k] M_n (x)$ is equal to the number of 
decorated permutations on $[n]$ which have $k$ weak excedences, i.e.
$[x^k] M_n(x)= A_{k,n}(1)$.  It would be interesting to find some 
$q$-analog of the above matrix whose permanent encodes 
$A_{k,n}(q)$.




\end{document}